# Intercultural Science-Art Project 'Magritte Meet Science'

Danilo Gregorin Afonso[a], Mary Oluwabunmi Akinade[b], Đorđe Baralić[c], Demian Nahuel Goos [d]*, Charles Gray[e], Abdulganiyu Jimoh[f], Lateef Olakunle Jolaoso[g], Jie Li[h], Sander Mack-Crane[i], Aung Zaw Myint[j], Kayode Oshinubi[k], Matheus Pires Cardoso[l], Dragana Radojičić[m]

[a]*University of São Paulo, São Paulo, Brazil;* [b]*Federal University of Technology, Akure, Nairobi, Nigeria;* [c]*Mathematical Institute SANU, Belgrad, Serbia;* [d]**Universidad Nacional de Rosario, Rosario, Argentina, demian@fceia.unr.edu.ar;* [e]*La Trobe University, Melbourne, Australia;* [f]*Adekunle Ajasin University, Akungba, Nigeria;* [g]*Sefako Makgatho Health Sciences University, Pretoria, South Africa;* [h]*Centrum Wiskunde & Informatica, Amsterdam, Netherlands;* [i]*University of California, Berkley, USA;* [j]*University of Mandalay, Mandalay, Myanmar;* [k]*Lagos State University, Lagos, Nigeria;* [l]*Universidade Federal do Tocantins, Tocantins, Brazil;* [m]*TU Wien, Vienna, Austria*

# Intercultural Science-Art Project 'Magritte Meet Science'

In this manuscript a joint project between young researchers from all around the world, the *Intercultural Science-Art Project*, is presented, in which all participants interpreted their own research with art. Art is meant in the most inclusive way possible, including drawings and digital works, but also songs and confectionary art. This project has two main goals. On the one hand it creates tools for young scientists to communicate their research to society by artistic means. On the other hand it tries to bring together different cultures and to strengthen ties between them.

## Introduction

We are a group of young scientists from all around the world who met in the Heidelberg Laureate Forum (HLF) 2019, a networking event which brings together young researchers in mathematics and computer science and laureates of these disciplines. As attendees of this forum, we were invited by Demian Nahuel Goos, a former HLF Alumnus who participated as blogger for the HLF 2019, to join his 'Intercultural Science-Art Project' and create artworks inspirend by our own scientific research.

The main goal of this project is to create awareness in society about the work of early-career researchers in academia and to make their work more visible. It is an indisputable fact that PhD theses—even though the peak of years of hard work and research—are habitually not read by anyone but the candidates, their advisors, and the committee. For enjoying such a read usually requires a basic comprehension of the theory that underlies these results or working at least on close and similar fronts. Either way, it is a small circle of scientists who have the privilege of fulfilling these requisites and people who do not work on academia can barely get an introductory notion of the content of the theses. It does not help this goal either that particularly in the presumably abstract fields of

Mathematics and Computer Science there are many prejudgments and misconceptions about their understandability and applicability to the real world.

So to make these works more tangible, it seemed natural to create a visual version of the theses and their main ideas, with which people could have a more concrete concept of what young researchers actually are working on and what the motivation behind their research is. This way the importance of their work would be easier to understand. Innovation and creativity are the key character traits for this activity since we expressed our research in an artistic fashion. Art is an interesting tool for the purpose of breaking the gap between science and society since it gives space to emotions and feelings, which again create a stronger bond between the scientific artist and the observer. The whole internal process of developing an artwork helped us young researchers in our own work as well, since we had to ask ourselves key questions about it, like "How can I explain my research comprehensibly?" We frequently do not ask these questions in everyday life, since we usually discuss our work with peers only, but they certainly help us to better understand what we are doing and what the goal of our work is.

A further element of the project was the intercultural focus. Keeping in mind that young researchers from all around the world meet at the HLF, it was a unique opportunity to see how people from different cultures interpret the task in their own way and how their cultural background, their country and language influence their way of facing this task. To highlight this, we wrote the titles of our artworks in our own language. This way all participants would get in touch with different languages and cultures and would create stronger bonds with people from all around the globe.

In this article we present the results of our project. After a brief presentation of the Leitmotiv, 'Magritte Meets Mathematics', one chapter is dedicated to each artwork with a brief explanation of the mathematical background. At the end, some final thoughts and commentaries are displayed. After some discussions limited to the HLF community (Goos, 2019), we hope to encourage futher debates and analysis with the mathematical community.

**The Leitmotiv 'Magritte Meets Science'**

René Magritte, a Belgian artist from the 20th century, is one of the most important representatives of surrealism, whose works in many cases can be given a scientific interpretation. Such is the case of *Le fils de l'homme*, in which we can see a man with a bowler hat—one of Magritte's key figures—whose face is covered by an apple. Magritte himself explained that the message behind this work was to describe the human desire to see beyond what can be seen. The problem is that, once we see the face behind the apple, it will not be enough: we will want to know what can be found behind the face. This is precisely the spirit of a scientist: Every obtained result gives place to new questions and new horizons and new theories. There will never be a last theorem. Behind every theorem, we can find another one.

However, we find in Magritte's work other ideas which are important to mathematics, like language. In this context, he addresses philosophical thoughts in a series of paintings in which we can see everyday objects, like hats and pipes, with different labels. The combination of both, the objects, and their labels, are always confusing, puzzling and irritating for the observer. The reason for this is that they initially seem to be contradictory. They challenge our intuition, our experience. However, there is always a clear idea behind this. The most famous one, *La Trahison des Images*, shows a pipe together with the text *Ceci n'est pas une pipe. —This is not a pipe*. With this work, Magritte discusses the difference between an object in the real world, its denomination, and its visual representation. Magritte asks what reality actually is and how we can describe our reality with language and images. He concludes that there are many limitations in language to represent reality.

Contrary to initial belief, science possesses such limitations as well. There are many paradoxes and counterintuitive results in mathematics and in this context, we may also wonder what kind of limitations mathematics possesses to faithfully describe our world. In addition to this, there are many misconceptions when it comes to society's understanding of science. Playing with these misconceptions society possesses about our work, Magritte's work turned out to be particulary interesting for this project. Magritte's inclusion of language in art helped us cover the intercultural approach as well.

**The artworks**

In this section, each artist presents her or his artwork with a brief statement about the scientific background. The titles are written in the young researcher's mother tongue, which is then adressed in brackets together with a translation into English.

**'Resolvendo sem resolver' (Portuguese, 'Solving without solving')**

The artwork (Figure 1) depicts a jigsaw puzzle with one piece missing. In a sense, the puzzle is solved, because the position of each piece is known, although it is not solved in a classical sense. Here, 'solution' entails to a general behavior of the missing piece, but not so much to its particularities.

The same ideas appear when one talks about the concept of weak solution to an elliptic partial differential equation.

Of course, given one such equation, a classical solution will be a function smooth enough that satisfies the equation in every point of the domain. However, differential equations are actually hard to solve. We can then forget about solutions in a classical sense and think about functions that satisfy the equation in some weaker sense. This makes the class of candidates to being a solutionmuch broader, hence it is easier to find existence theorems. Now, the beautiful thing is that under some conditions these weak solutions will be solutions in the classical sense. Such conditions are studied by regularity theory, see Badiale & Serra (2011), Brezis (2010) or Kesavan (1989).

We thus have solved the problem, in some non-orthodox sense. Just like in the artwork.

**'Nkan tio ba kọjaa mathematics, o kọjaa eniyan' (Yoruba, 'Whatever goes beyond mathematics, goes beyond man')**

Looking externally at the work (Figure 2), one may be quick to conclude that it's the map of Nigeria, my country of origin. This image, however, does not include the two rivers which divides the country into three regions, hence visualizing the country as a single unit.

The two lines contained in the image are representations of few of the numerical simulation results performed on a Lassa fever model which I developed in the course of my masters' degree research. Lassa Fever is a highly fatal acute Viral Hemorrhagic fever caused by Lassa Virus. The first occurrence of this disease was in Nigeria in 1969 and presently, Nigeria is home to over 50% of the total human Lassa fever population.

The yellow line in the image represents the numerical simulation result of the model's infectious class without the incorporation of any control strategy while the red line represents the numerical simulation result of the model's infectious class with the incorporation of three different control strategies.

The phrase "Nkan tio ba kọjaa mathematics, o kọjaa eniyan" is written in Yoruba, my native language. The English translation is "Whatever goes beyond mathematics, goes beyond man". This is to simply state that Mathematics is a tool which helps to proffer solutions to all human problems, provided its properly maximized.

**'Квазиторусна 4-многострукост Quasotoric 4-manifold' (Serbian, 'Quasitoric manifolds')**

Quasitoric manifolds and small covers, their real analogues appeared in a seminal paper (Davis & Januszkiewicz, 1991) of Davis and Januszkiewicz as a topological generalizations of non-singular projective toric varieties and real toric varieties. The manifolds have a locally standard $(S^d)^n$ action where $d = 0$ in the case of small covers and $d = 1$ in the case of quasitoric manifolds, such that the orbit space of the action is identified with a simple polytope as a manifold with corners. The simplest examples are manifolds over the *n*-dimensional simplex $\Delta^n$, $\mathbb{C}P^n$ for quasitoric manifolds and $\mathbb{R}P^n$ for small covers and they are unique up to homeomorphism.

In the last decades, toric topology experienced an impressive progress. The most significant results are summarized in recent remarkable monograph (Buchstaber & Panov, 2015). However, one of the most interesting problems in toric topology such as classification of simple polytopes that can appear as the orbit spaces of some quasitoric manifolds and classification of quasitoric manifolds and small covers over a given simple polytope are still open. In my PhD thesis, I studied embeddings of quasitoric manifolds into Euclidean spaces (Baralic, 2014) and the mapping degrees of maps between 4 dimensional quasitoric manifolds (Baralic & Grujic, 2016).

Let $P^n$ be a simple polytope with $m$ facets $F_1,\ldots,F_m$. By Definition of quasitoric manifold (Buchstaber & Panov, 2015) it follows that every point in $\pi^{-1}(rel.int(F_i))$ has the same isotropy group which is a one-dimensional subgroup of $(S^1)^n = T^n$ We denote it by $T(F_i)$. Therefore, for a point $p \in rel.int(F_i)$ we have that $\pi^{-1}(p)$ is homeomorphic to a torus $T^{n-1}$. Similarly, for a point $q$ lying in the interior of some $k$ faces we have that $\pi^{-1}(q)$ is homeomorphic to a torus $T^{n-k}$ and the fixed points of $T^n$ action on quasitoric manifold $M^{2n}$ correspond to the vertices of $P^n$.

This description of the orbits of $T^n$ action on a quasitoric manifold may give us some insights about topology and geometry of quasitoric manifolds and an idea how to visualize them. The artwork (Figure 3) was inspired by Figure 4.5 of my Ph.D. Thesis (Baralic, 2013) and acquires this concept for $n = 2$ and a 4 dimensional quasitoric manifolds, where the orbit space is a polygon and a torus $T^2$ is the pre-image of any point inside the polygon and a circle $S^1$ is pre-image of any interior point of an edge of the polygon.

**'Hic est draco' (Latin, 'Here is a dragon')**

The P vs. NP conjecture in complexity theory (Cook, 1971) is one of the most prominent open problems in mathematics and computer science. There are many reasons for this. It is reasonably easy to get an intuitive idea of what it states, it has many deep implications to everyday life in both, science and society, and last but not least, it was included by the *Clay Mathematics Institute* in the selected group of *Millenium problems*.

An artistic approach to this problem should faithfully visualize what it represents to the mathematical community. What we feel when we think about it and what a solution to the conjecture would mean. If the known mathematics is a map, then conjectures like this one represent the limits of the known world. Something human mankind wishes to know and understand better. We ignore what can be found behind these areas and the mere idea of facing them seems to be a perilous and foolish endeavor. The first world maps labeled unexplored territories with the inscription *Hic sunt dracones – Here be dragons*, together with the depiction of sea monsters and mythological creatures. Now, what is an open conjecture if not a mysterious, untamed monster waiting to be conquered? We might think that we understand it, but truth is that only the most skilled adventurers may have a serious chance against this beast, and the glory that comes with a definitive answer to the P vs. NP issue is comparable to the glory of overcoming a beast in tall tales.

In classical medieval style, my monster (Figure 4) is a mix of different real-world animals. This can frequently be observed in tapestry all around Europe, where the artist creating the tapestry depicted strange and exotic landscapes without ever visiting them. Basing their art on descriptions of explorers, strange creatures emerged from this interaction. The inscription in Spanish, '*This monster cannot be tamed in polynomial time*'.

**'incorrigible tidy::vert %>% exegesis()' (Programming language *R*)**

In his 1929 surrealist painting, *The Treachery of Images*, René Magritte declares *Ceci n'est pas une pipe* (This is not a pipe). In so doing, he highlights this is but an image, a representation, of a pipe, not truly a pipe itself. I recently wrote a song (Gray, 2019) with lots of references to R and the tidyverse:: metapackage (Wickham,. 2017). Here, I will unpack the pipe operator, %>%, that features throughout the lyrics (Figure 5). The %>% pipe operator is the first major concept introduced in Wickham & Grolemund (2016).

What is an operator? We often forget that operators are, themselves, functions. For example, + is a function that takes two arguments, numbers, and returns a single number. Algebraically, $3 + 2 = 5$ is shorthand for $+(3,2) = 5$. For those with formal mathematical training, multiple uses of the %>% operator in a single line of code can be thought of in terms of a coding instantiation of a composite of functions. What is a composite? Let f and g be real functions. The composite of $f$ with $g$ is the real function $g \circ f$ given by the formula $(g \circ f)(x) \coloneqq g(f(x))$. For reasons that only made sense to me once I reached graduate-level mathematics, we read a composite of functions from right to left. And just to break our brains a little, algebraically, the composite operator is a function, so we have $g \circ f = \circ (f, g)$! The pipe, %>%, operator is the R-language equivalent to the composite $\circ$ operator on real functions.

Why do I love to %>%? Here is an example with three functions: $(h \circ g \circ f)(x) \coloneqq h(g(f(x)))$.

```
set.seed(39)
    # get a random sample size between 20 & 100
    sample(seq(20, 100), 1) %>% # this f(x) goes into
      # generate sample from normal distribution with
      # mean 50 & sd 0.5
    rnorm(., 50, 0.5) %>% # g, so, now g(f(x)), which goes into
      # calculate mean of that sample
      mean() # h, so h(g(f(x)))
```

To see how this is the $(h \circ g \circ f)(x)$ instantiation, reading from right to left, we take a look at the $h(g(f(x)))$ instantiation of the same code.

```
    # this line of code is equivalent to above
    # h(g(f(x))) is less text
    # but the algorithm is harder to ascertain
    mean(rnorm(sample(seq(20, 100), 1), 50, 0.5))
```

The reader is invited to consider if they agree with the author that it is harder to read the symbols so close together, in this $h(g(f(x)))$ instantiation of the code. Also, arguably more importantly, one does not have the ability to comment each component of the algorithm. There is a downside to the %>%, however. The longer a composite becomes, the more difficult it is to identify errors.

**'Eyi kii ṣe ọlọje x-ray awọn pixels.' (Yoruba, 'This is not x-ray scan, it's pixels')**

During data collection in my research work I noticed that at every digital image taking section in radiography unit in hospitals, radiologist doctor's has the strong believe that they are truly diagnosing patients for every image scan to detect illness, but this doctor's aren't realising that what they think they are diagnosing are series of pixels representing or depicting the real medical image of patients. Therefore, my artwork (Figure 6) reads *'Eyi kii ṣe ọlọje x-ray awọn pixels.', 'This is not x-ray scan it's pixels.'*

**'Èyi kọ nṣe pọinti aiyipada' (Yoruba, 'This is not a fixed-point theorem')**

In my artistic photograph (Figure 7), I depict the field of fixed-point theory, which has become one of the most dynamic areas of research in the last five decades (Brezis, 2010). It is a beautiful mixture of analysis (pure and applied), topology and geometry and has played a basic role in the development of several theoretical and applied fields such as nonlinear analysis, integral and differential equations, dynamical system, optimization, physics, economics and engineering.

Several top-rated international journals have been devoted for publications on fixed-point theory, this include Journal of fixed point theory and application, Fixed-point theory and application, Fixed point theory, Advances in fixed point theory and, Recent advances in fixed point theory and applications, and new developments in fixed point-theory and applications.

The field has further played a key role in proving the existence and uniqueness of solutions of differential and integral equations. It also presents the conditions under which single-valued and multivalued mappings have solutions. A celebrated bedrock on this field is the work of Stefan Banach in 1922, and since then, several fixed-point results have been announced by many authors.

### '三十三' (Chinese, 'Thirty-three')

This (Figure 8) is a matcha (抹茶, Japanese green tea) chocolate pie, decorated with 'forest-like' matcha cream. The person who requested me to make this pie, would like to have a 'magic square' on it. One of the first magic squares depicted in art was the 4×4 square in an engraving entitled *Melancholia I* by Albrecht Dürer. Filled with numbers 1 through 16, the Dürer Square demonstrates the magic constant, 34 and the year of engraving, 1514, in the middle of the bottom row of the square.

This magic square on the pie is the same one as on the *Passion façade* of Antoni Gaudí's *Sagrada Família*. It also has a religious meaning as 'Jesus Crucifixion Age 33'. Although it is not a true magic square because the grid contains the duplicate numbers of 10 and 14, it does consist of a magic constant--The number 33, with sums of rows, columns, diagonals, and 2 x 2 sub-squares all equal to 33. The reason why this person wanted to have this special pie is that it was customized for his girlfriend, whose name sounds like '33' in Chinese.

### 'This is not two worlds. It's one.'

Number theory is marked by astonshing connections between seemingly unrelated ideas. To understand how common prime numbers are, we translate the question into the entirely different language of complex analysis. To prove Fermat's last theorem—a simple statement about the integers—we need to combine powerful tools in complex analysis as well as algebraic geometry. The Langlands program is a modern example of such an astonishing connection (which in fact generalizes many of the ideas in the proof of Fermat's last theorem). One view of the Langlands program is that certain geometric spaces called Shimura varieties should have a description in terms of analytic objects called automorphic representations. In order to avoid too much technical detail, we can simply think of this as a deep connection between two different worlds: the geometric world and the spectral world (which is our name for the world of analysis, the word "spectral" distantly related to the "spectrum" of light). In the artwork (Figure 9), we see two very different and very separate worlds—our spectral and geometric worlds. Decades of work by many people has led to the development of a fruitful approach to establishing this connection, known as the Langlands-Kottwitz method. A prime example is the paper Kottwitz (1992). This method relies on an extremely complicated machine called the trace formula. The trace formula is an equation which straddles the divide: the two sides of this equation are called the spectral side and the geometric side, and they express the same information in the languages of these two different worlds. In the artwork we see our two worlds are connected by a tremendous and complicated machine—the trace formula, which itself has a spectral side and geometric side. To get information from one world to the other, we first have to establish a connection between the geometric world and the geometric side of the trace formula; then the trace formula can be used to transfer this information to its spectral side, and finally the spectral side of the trace formula is connected to the spectral world.

My own research is concerned largely with the first step in this process, describing geometric spaces in a way the geometric side of the trace formula can understand—laying cables from the geometric world to the geometric side of the trace formula. This process is often referred to as "point-counting", because the main step is to get a good description of all the points of our space. The papers Kisin (2017) and Shin (2009) carry out this process for certain classes of geometric spaces, and my work in progress combines the ideas from these papers to handle a new class of spaces—building on the infrastructure they have established to expand the network of connections ever wider.

### 'Predator-prey Model ရဲ့အရေအတွက် ပေါက်ကွဲခြင်း စနစ်' (Burmese, 'Population Dynamics of Predator-prey System')

This Pheasant Math-Art (Figure 10) is referred to the predator-prey relationship system graphs of the nonlinear of differential equations in mathematical biology or ecology.

In the back part of Pheasant, presented result graph between the lynx and the snowshoe hare the predator prey relationship. The yellow shows the population of lynx, while the red shows the population of hares. At the start of the graph, the lynx population was very high, which the hare population was relatively low (i.e. the snowshoe hare forms a large staple in the lynx diet. Without the hare, the lynx would starve. However, as the lynx eats the hare, or many hares, it can reproduce. Thus, the lynx population expands. With more lynx hunting, the hare population rapidly declines.)

For Lotka-Volterra equations (Murray, 2002 and 2003; Britton, 2012), the Pheasant tails are referred to critical point in the first quadrant could be an asymptotically stable node or spiral point when the limiting population of the prey species in the absence of the predator species.

The front part of Pheasant referred Prey-Predator dynamics as described by the level curves of a conserved quantity for Lotka-Volterra equations.

In the neck of Pheasant said that Kermack-McKendrick model of propagation of infectious disease. A short calculation shows that $x(t)$ converges to a constant, say $x(t) \to x^*$ where $x^*$ can be found by solving the equation $C = rx^* - d\ln x^*$.

### 'Kiraki' (Yoruba, 'Cracks')

The artwork (Figure 11) is all about cracks in structures as we know that cracks leads to collapse of structures. The idea is to proffer solutions to cracks in structures in order to predict the location, size and depth of a crack.

**'A caixa de Cantor' (Portugese, 'Cantor's Box')**

Paradoxes are fascinating objects, especially when they challenge our intuition and instigate us to reflex and makes us comprehend what plays with our perception and understanding. The art has been a fellow of paradoxes, it is in this way that the paradoxes sentences embody in lines, traces and colours giving another vison of the contradictions, intuitions errors and situations which produce paradoxes for us finite beings showing the other face of the same coin. It is in this purpose that the artwork (Figure 12) *Caixa de Cantor* tries to show another facette of a famous paradox from set theory, the wood box floating in a blue surface talks about infinite sets and the transfinite cardinals.

First of all, we would like to describe the process of the production of this artwork. Inspired in the artwork of René Magritte, Cantor's Box was formed through the free app Autodesk Sketchbook which has a notable amount of drawing tools for beginning artists. The sketch and the final result was carefully made during a couple of hours through the touchscreen of a smartphone, giving a body, even in a very simplistic way, to the *Hilbert's Hotel Paradox* (Stillwell, 2010; Figueiredo, 2011), beyond the box there is a short phrase written in Portuguese, the mother tongue of the artist, which assists as hint to remember the paradox that happens in the transfnite arithmetic, beyond that there is a precise suggestion given by the title of the artwork.

So, what is Cantor's Box about? Well, let us make some clarifications, assume that our box is the set of the natural numbers, i. e. any positive and integer number that we can think is in our box, and let us imagine that each number is a block in a cube form, and all the infinite number are organized in piles so that our box is 'full'. Therefore, the block 1024 is there inside our box as like as the block 123456 is as well. Then, we can make some conclusions about our box, first, it is infinite and second, all the naturals have a place inside it, and probably, there is no place for more. But Cantor would disagree with that: suppose we receive the block 0 and we want to save it in the box, we probably think that it is impracticable because our box is full, but Hilbert has a solution for that. He will advise us to pick the block 1 and put it the 2's place, and block 2 in the 3's place and so on. Arranging in this way we will have an additional place for our block 0, the place which was occupied by block 1. It appears simple, even though our box is 'full', it allows room for one block more.

But now, suppose that have an infinite collection of numbers blocks labeled with even numbers and we want to save them as well. Now, we should organize our box so that the block of number $x$ moves to the block of number $f(x) = 2x + 1$. In this way, we will have all the places which were occupied by the even blocks empty, and we can put the new blocks inside our box. That is possible because our function $f(x) = 2x + 1$ is a one-to-one correspondece. That justifies the phrase *'it is full, but there is space'*.

**'Knjiga limitiranih naloga' (Serbian, 'The Limit Order Book')**

As the sun goes up, we wake up wanting to know in which direction in the forest we should go in order to harvest some dollars from the tree. Should I plant, pick, or just lay on the grass to get some tan? Should I buy, sell, or idle? Stock market and high-frequency trading have attracted attention from the financial institutions and hedge funds, but also from academics, for example see Radojičić & Kredatus (2020). Before the sun goes down, we should decide which signs and indicators we take to trust them and believe in them.

The text in the picture (Figure 13) is from the song: 'Tesla' in the performance of the group Letu Štuke, and it means: 'because the world is in ourselves, not in the towers'.

As Einstein once said: 'You can't blame gravity for falling in love'. We can say that you should not blame mathematics for not doing well on stock markets.

**Conclusions and final thoughts**

The Intercultural Science-Art Project was an enriching experience for all participants. During the HLF we had several opportunities, among those a Science-Art session organized by the HLF team, to discuss the artworks and our thoughts:

- The thinking process during the creation of our artworks was similar to the one we use on an everyday basis in our research. You ask questions like "How can I solve this?" or "Is this really consistent?" After all, in both cases, we are talking about a creative instance which may even need a moment of inspiration and a carefully approached and well-thought rational instance. This connection between arts and science can be used to demystify the work of a scientist.

- Science and art are actually quite similar. They both take our reality—our real world—and interpret it by creating an abstraction of it. While one is heavily focused on objective preciseness and a clear and uniformly accepted scientific method, the other one is a rather subjective, emotional interpretation of the same affairs. The interesting thing about this is that this means we can make science more attractive to more people since it implies that scientists are nothing more than creative artists who just happen to use their own, particular tools and methods. But then, what artist does not?

- Mathematics and Computer Sciences are described as really abstract fields, whereas arts are considered to be quite visual. However, arts can be really abstract as well, whereas science is heavily inspired by the real world we live in. So, this leads to the mistake of thinking that they both work in opposite directions, while both—arts and science—have one common motivation: the search for beauty and harmony.

- In both fields, science and arts, the relationship between the scientist/artist and user/observer is very similar. The first presents a creation of his own, while the latter uses, interacts and interprets it. Both, the scientists and the artists have a very clear intention for their creation, but it is up to the user and the observer to decide what to do with it. Of course, the more interpretations a creation can be given and the more applications it has, the more valuable this creation will be for society.

- One important aspect of intertwining mathematics with arts is that it is a way to make mathematics more visible to the non-scientific community due to the strong bond between arts and human emotions. This way people can get in touch with science in a more natural way. By reaching people's hearts and focusing on their critical thinking at the same time, we create a strong bond between apparently unconnected elements of human nature.

On a final note, both—science and arts—have one thing in common. If not corrupted, they are egalitarian, democratic, and peaceful pillars of society and its development. They both represent what humanity should be striving for.

**Acknowledgements**

We would like to thank the whole HLF team for their unconditional support and enthusiasm concerning this project. In particular, we would like to thank the Young Researchers Relations and the Communications & Media group. This project would not have been possible without their work.

**Figures**

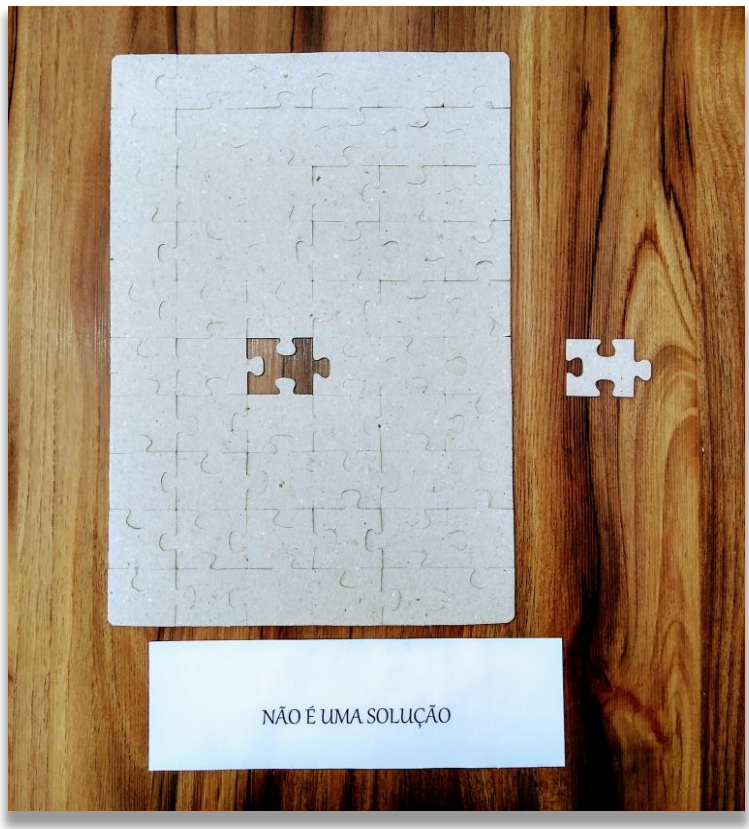

Figure 1. 'Resolvendo sem resolver', (Portuguese, 'Solving without solving'), Danilo Gregorin Afonso, 2019, Jigsaw puzzle.

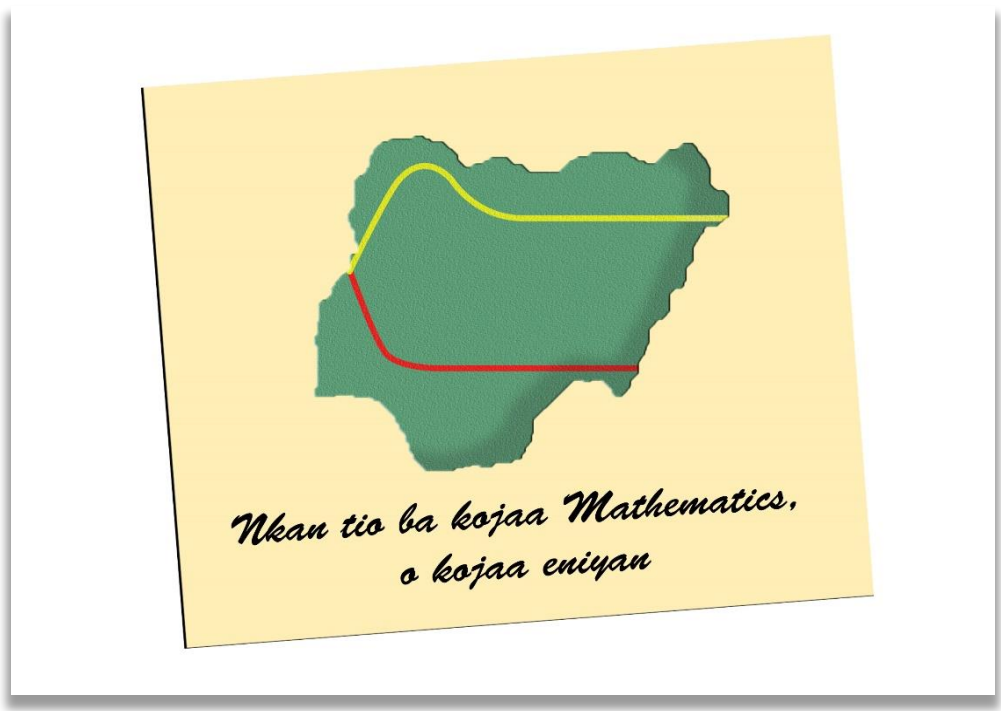

Figure 2. 'Nkan tio ba kọjaa mathematics, o kọjaa eniyan' (Yoruba, 'Whatever goes beyond mathematics, goes beyond man'), Akinade Mary Oluwabunmi, 2019, digital.

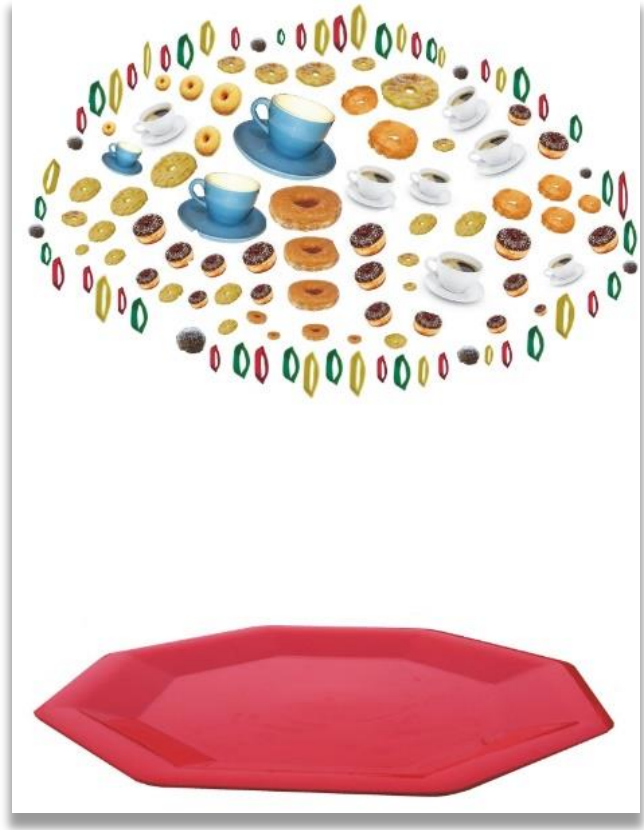

Figure 3. 'Квазиторусна 4-многострукост Quasotoric 4-manifold' (Serbian, 'Quasitoric manifolds'), Đorđe Baralić, 2019, digital collage.

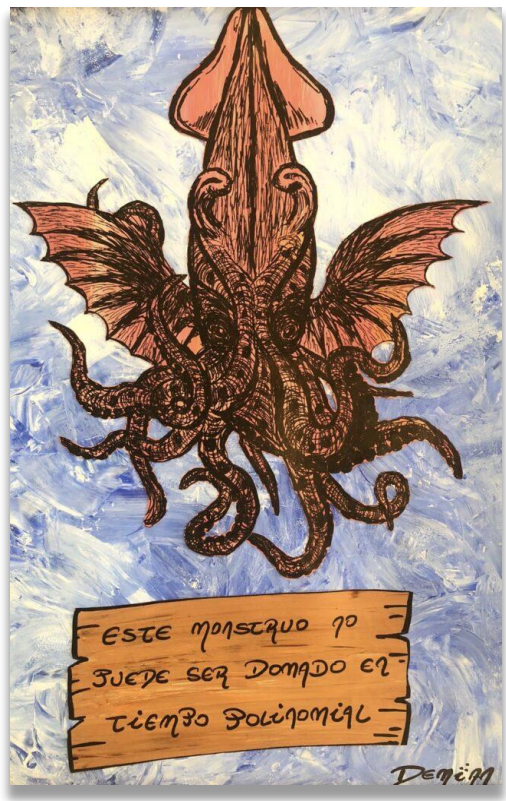

Figure 4. 'Hic est Draco' (Latin, 'Here is a dragon'), Demian Nahuel Goos, 2019, ink and acrylic on paperboard.

Figure 5. Lyrics of 'incorrigible tidy::vert %>% exegesis()', Charles Gray, 2019, song.

Figure 6. 'Eyi kii șe ọlọjẹ x-ray awọn pixels.' (Yoruba, 'This is not x-ray scan, it's pixels.'), Jimoh Abdulganiyu, 2019, drawing.

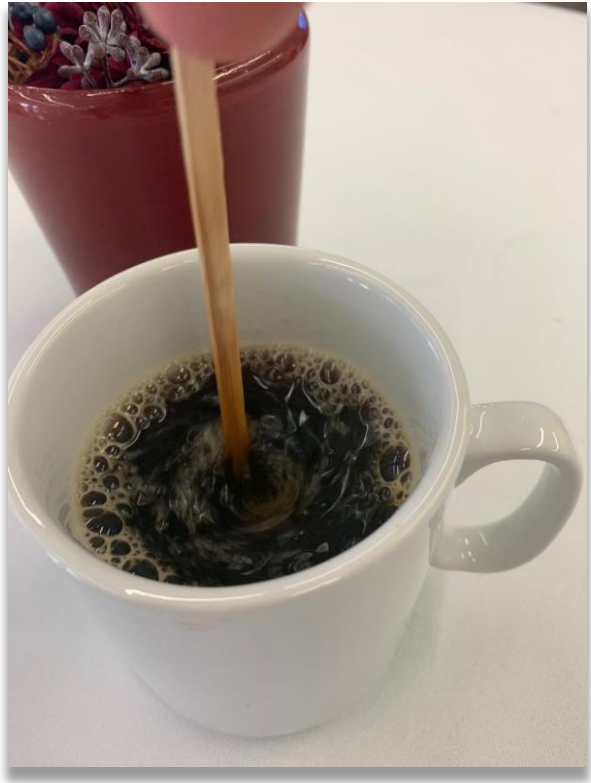

Figure 7. Èyi kọ nṣe pọinti aiyipada' (Yoruba, 'This is not a fixed-point theorem'), Lateef Olakunle Jolaoso, 2019, photography.

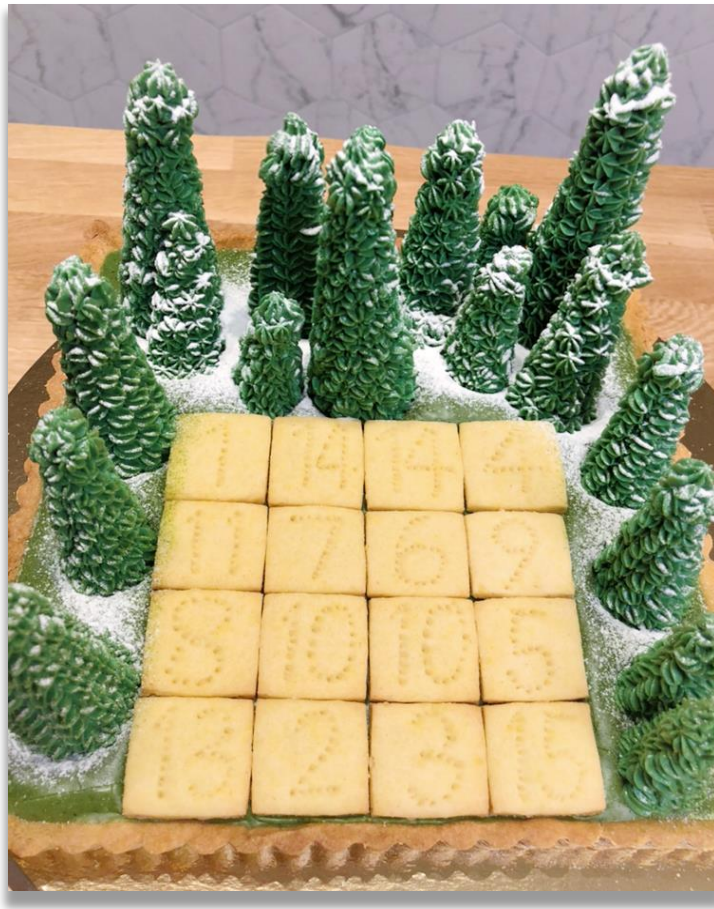

Figure 8. '三十三' (Chinese, 'Thirty-three'), Jie Li, 2019, confectionary art.

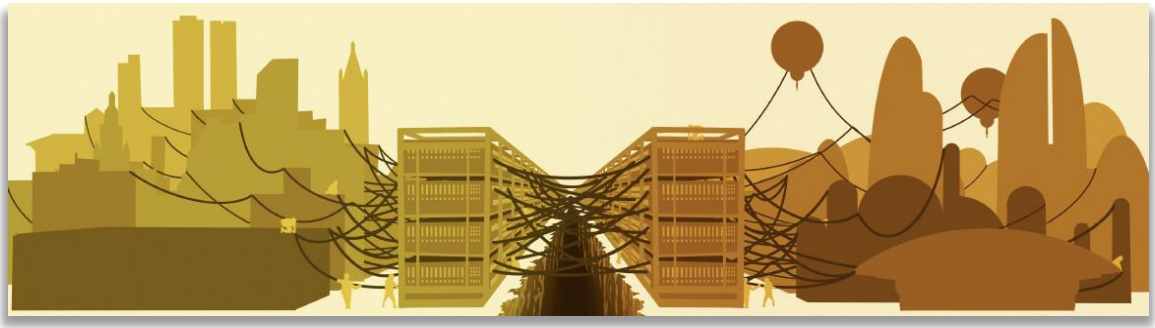

Figure 9. 'This is not two worlds. It's one.', Sander Mack-Crane, 2019, digital.

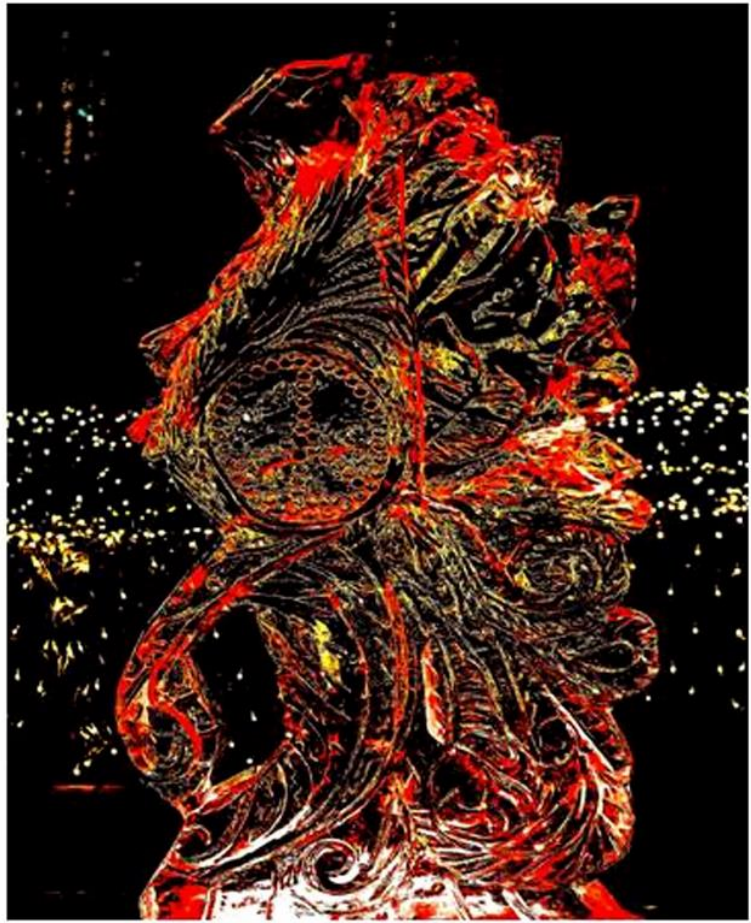

Figure 10. 'Predator-prey Model ရဲ့ အရေအတွက် ပေါက်ကွဲခြင်း စနစ်' (Burmese, 'Population Dynamics of Predator-prey System'), Aung Zaw Myint, 2019, digital.

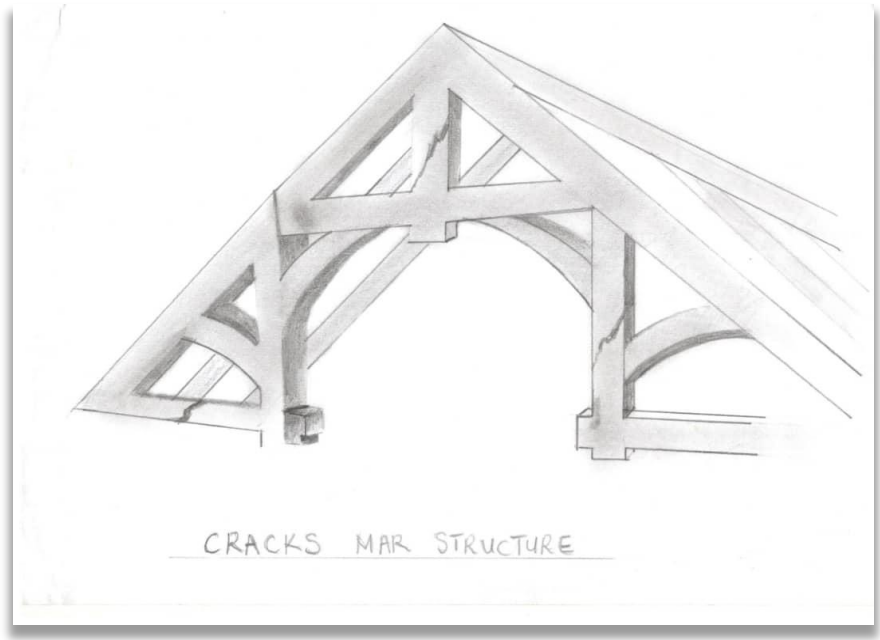

Figure 11. 'Kiraki' (Yoruba,'Crack'), Oshinubi Kayode, 2019, drawing.

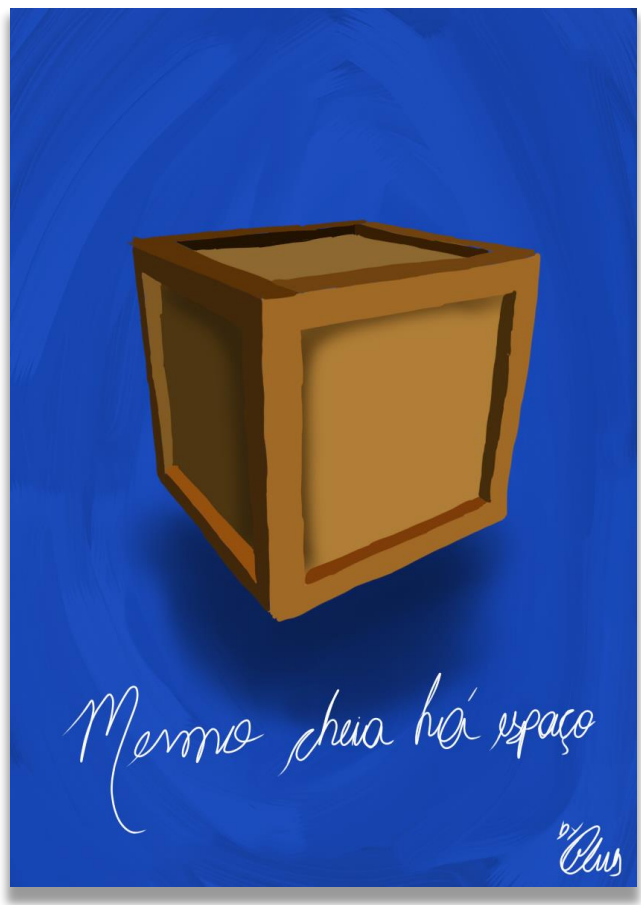

Figure 12. 'A caixa de Cantor' (Portugese, 'Cantor's box'), Matheus Pires Cardoso, 2019, digital.

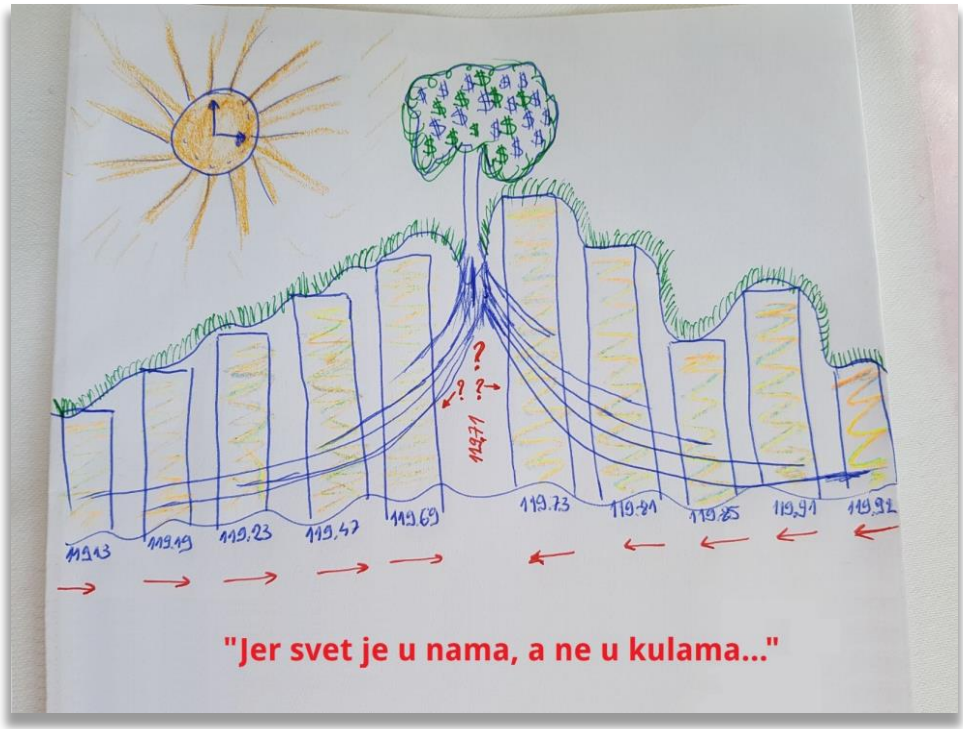

Figure 13. 'Knjiga limitiranih naloga' (Serbian, 'The Limit Order Book'), Dragana Radojičić, 2019, pen on paper.